\documentstyle[11pt]{article}

\newtheorem{thm}{Theorem}[section]
\newtheorem{rk}[thm]{Remark}
\newtheorem{prop}[thm]{Proposition}
\newtheorem{clly}[thm]{Corollary}
\newtheorem{lemma}[thm]{Lemma}

\begin{document}

{\centerline{\LARGE{\bf \underline{Generalized fixed-point algebras }}}
{\centerline{\LARGE{\bf \underline{of certain actions on crossed products.}}}
\vspace{.2in}
\centerline{\large{Beatriz Abadie\footnote{This work is a generalization of
part of the contents of the author's Ph.D dissertation submitted to the
University of California at Berkeley in May 1992.
\newline 1991 MR Classification: Primary 46L55; Secondary 46L80, 46L87.}}}
\vspace{.2in}
{\small{\em{\underline{Abstract:} Let $G$ and $H$ be two locally compact groups
acting on a C*-algebra $A$ by commuting actions $\lambda$ and $\sigma$. We
construct an action on $A\times_{\lambda}G$ out of $\sigma$ and a unitary
2-cocycle $u$. For $A$ commutative, and free and proper actions $\lambda$ and
$\sigma$, we show that if the roles of $\lambda$ and $\sigma$ are reversed, and
$u$ is replaced by $u^{*}$, then the corresponding generalized fixed-point
algebras, in the sense of Rieffel, are strong-Morita equivalent. We apply this
result to the computation of the K-theory of quantum Heisenberg manifolds.}}}
\vspace{.3in}
\newline\underline{\large{\bf Introduction.}}
 Given two commuting actions $\lambda$ and $\sigma$ of locally compact groups
$G$ and $H$, respectively, on a C*-algebra $A$, we study the action
$\gamma^{\sigma ,u}$ of $H$ on $A\times_{\lambda}G$ defined by
\[(\gamma_{h}^{\sigma ,u}\Phi)(x)=u(x,h)\sigma_{h}(\Phi(x)),\]
where $\Phi\in C_{c}(G,A)$, $h\in H$,  $x\in G$, $u(x,h)$ is a unitary element
of the center of the multiplier algebra of $A$, and $u$ satisfies the cocycle
conditions
\[u(x_{1}x_{2},h)=u(x_{1},h)\lambda_{x_{1}}(u(x_{2},h))\mbox{\hspace{.1in}and\hspace{.1in}}u(x,h_{1}h_{2})=u(x,h_{1})\sigma_{h_{1}}(u(x,h_{2})).\]
The study of this situation was originally motivated by the example of quantum
Heisenberg manifolds (\cite{rfhm}), which can be described as generalized
fixed-point algebras (\cite{rfpa}) of actions of this form, when
$A=C_{0}(R\times T)$, and $G=H=Z$.

This work is organized as follows. In Section 1 we define the action
$\gamma^{\sigma ,u}$ and show that for $G$ and $H$ second countable, and $A$
separable, the crossed product $A\times_{\lambda}G\times_{\gamma^{\sigma ,u}}H$
is isomorphic to a certain twisted crossed product of the algebra $A$ by the
group $G\times H$.

In Section 2 we assume that the algebra $A$ is commutative and show that for
free and proper actions $\lambda$ and $\sigma$, the generalized fixed-point
algebra of $A\times_{\lambda}G$ under $\gamma^{\sigma ,u}$ and that of
$A\times_{\sigma} H$ under $\gamma^{\lambda ,u^{*}}$ are strong-Morita
equivalent.

In Section 3 we apply these results to show that the K-groups of the quantum
Heisenberg manifolds do not depend on the deformation constant. This enables us
to compute them, by calculating them in the commutative case.

In what follows, for a C*-algebra $A$, ${\cal M}(A)$ denotes its multiplier
algebra, ${\cal Z}(A)$ its center and ${\cal U}(A)$ the group of unitary
elements in $A$.   All actions of locally compact groups on C*-algebras are
assumed to be strongly continuous.
All integrations on a group $G$ are with respect to a fixed left Haar measure
$\mu_{G}$ with modular function $\Delta_{G}$.
\section{Actions on crossed products.}
\label{gen}
For locally compact groups $G$ and $H$ acting on a C*-algebra $A$ by commuting
actions $\lambda$ and $\sigma$, respectively, and a 2-cocycle on $G\times H$,
we define an action $\gamma^{\sigma ,u}$ of $H$ on $A\times_{\lambda}G$. We
show in Proposition \ref{tcp} that, when $A$ is separable, and $G$ and $H$ are
second-countable, the crossed product $A\times_{\lambda}G\times_{\gamma^{\sigma
,u}}H$ is a twisted crossed product of $A$ by $G\times H$.
\begin{prop}
\label{aut}
Let $G$ be a group acting on a C*-algebra $A$ by an action $\lambda$, and let
$v:G\rightarrow{\cal U Z M}(A)$ verify the cocycle condition
\[ v(xy)=v(x)\lambda_{x}(v(y)).\]
Let $\sigma\in Aut(A)$ commute with $\lambda$, and, for $\Phi\in C_{c}(G,A)$,
define
\[(\gamma^{\sigma ,v}\Phi)(x)=v(x)\sigma(\Phi(x)).\]
Then $\gamma^{\sigma ,v}$ extends to an automorphism on $A\times_{\lambda}G$.
\end{prop}
\vspace{.1in}
\underline{\em{Proof:}}

\vspace{.2in}
Let $(\Pi,V)$ be a covariant representation of the C*-dynamical system
$C^{*}(G,A,\lambda)$ on a Hilbert space ${\cal H}$, and let $\Pi\times U$
denote its integrated form.
Let $\Pi^{\sigma}$ denote the representation of $A$ on ${\cal H}$ defined by
$\Pi^{\sigma}(a)=\Pi(\sigma(a)),$
and let $\tilde{V}$ be the unitary representation of $G$ on ${\cal H}$ given by
$\tilde{V}_{x}=\Pi(v(x))V_{x}$, where $\Pi$ also denotes its extension to
${\cal M}$.
Then $(\Pi^{\sigma},\tilde{V})$ is a covariant representation of
$C^{*}(G,A,\lambda)$: for $x\in G$, and $a\in A$ we have
\[\tilde{V}_{x}\Pi^{\sigma}(a)\tilde{V}_{x^{-1}}=\Pi(v(x))V_{x}\Pi(\sigma(a))\Pi(v(x^{-1}))V_{x^{-1}}=\]
\[=\Pi(v(x))\Pi(\lambda_{x}\sigma(a))V_{x}\Pi(v(x^{-1}))V_{x^{-1}}=\]
\[=\Pi(v(x))\Pi(\sigma\lambda_{x}(a))\Pi(\lambda_{x}v(x^{-1}))=\Pi^{\sigma}(\lambda_{x}(a)).\]
We now show that for $\Phi$ in $C_{c}(G,A)$ we have that $(\Pi\times
V)(\gamma^{\sigma,v}\Phi)=$
\newline $=(\Pi^{\sigma}\times \tilde{V})(\Phi)$, which ends the proof: for any
$\xi$ in ${\cal H}$, we have
\[[(\Pi\times
V)(\gamma^{\sigma
,v}\Phi)](\xi)=\int_{G}\Pi[(\gamma^{\sigma,v}\Phi)(x)]V_{x}\xi
dx=\int_{G}\Pi(v(x))\Pi[(\sigma\Phi)(x)]V_{x}\xi dx=\]
\[=\int_{G}\Pi^{\sigma}[\Phi(x)]\tilde{V}_{x}\xi dx=[(\Pi^{\sigma}\times
\tilde{V})(\Phi)](\xi).\]
\begin{flushright}
Q.E.D.
\end{flushright}
\vspace{.2in}
\begin{prop}
\label{act}
Assume that $G$, $\lambda$, and $A$ are as in Proposition \ref{aut} and that H
is  a locally compact group acting on $A$ by an action ${\sigma}$ commuting
with ${\lambda}$. Let
\[u:G\times H\rightarrow {\cal UZM}(A)\]
be continuous for the strict topology in ${\cal M}(A)$, and satisfy
\[u(xy,h)=u(x,h)\lambda_{x}u(y,h)\mbox{\hspace{.1in}and\hspace{.1in}}u(x,hg)=u(x,h)\sigma_{h}u(x,g),\]
for $x,y\in G$ and $h,g\in H$.
For $h\in H$ and $\Phi\in C_{c}(G,A)$, let
\[(\gamma^{\sigma ,u}_{h}\Phi)(x)=u(x,h)\sigma_{h}(\Phi(x)).\]
Then $h\mapsto \gamma_{h}$ is a (strongly continuous) action of $H$ on
$A\times_{\lambda}G$.
\end{prop}
\vspace{.1in}
\underline{\em{Proof:}}

\vspace{.2in}
By Proposition \ref{aut} we have that $\gamma_{h}\in Aut(A\times_{\lambda}G)$,
for all $h\in H$. Besides, the cocycle condition implies that $\gamma^{\sigma
,u}_{h_{1}h_{2}}\Phi(x)=\gamma^{\sigma ,u}_{h_{1}}\gamma^{\sigma
,u}_{h_{2}}\Phi(x)$.
\newline Finally, $h\mapsto \gamma^{\sigma ,u}_{h}\Phi$ is continuous for any
$\Phi\in C_{c}(G,A)$:
\[ \|\gamma^{\sigma ,u}_{h}\Phi-\gamma^{\sigma
,u}_{h_{0}}\Phi\|_{A\times_{\lambda}G}\leq\|\gamma^{\sigma
,u}_{h}\Phi-\gamma^{\sigma ,u}_{h_{0}}\Phi\|_{L^{1}(G,A)}=\]
\[=\int_{G}\|u(x,h)\sigma_{h}(\Phi(x))-u(x,h_{0})\sigma_{h_{0}}(\Phi(x))\|_{A}dx\leq\]
\[\leq
\int_{supp(\Phi)}\|\sigma_{h}(\Phi(x))-\sigma_{h_{0}}(\Phi(x))\|_{A}+\|(u(x,h)-u(x,h_{0}))\sigma_{h_{0}}(\Phi(x))\|_{A}dx,\]
which converges to $0$ when $h$ goes to $h_{0}$, because $u$ is continuous, and
$\sigma$ is strongly continuous.
\begin{flushright}
Q.E.D.
\end{flushright}
Next Proposition shows that the double crossed product
$A\times_{\lambda}G\times_{\gamma^{\sigma ,u}}H$ is isomorphic to a twisted
crossed product. Since twisted crossed products are defined for separable
algebras and second-countable groups, we add these conditions.
\vspace{.2in}
\begin{prop}
\label{tcp}
Let $G$, $H$, $A$, $u$, $\lambda$, $\sigma$ and $\gamma^{\sigma ,u}$ be as in
Proposition \ref{act}. If $A$ is separable and $H$ and $G$ are
second-countable, then $A\times_{\lambda}G\times_{\gamma^{\sigma ,u}}H$ is
isomorphic to the twisted crossed product $A\times_{(\lambda,\sigma),U}(G\times
H)$, where
\[(\lambda,\sigma)_{(x,h)}(a)=\lambda_{x}\sigma_{h}(a)\mbox{\hspace{.1in}and\hspace{.1in}}U((x_{0},h_{0}),(x_{1},h_{1}))=\lambda_{x_{0}}(u(x_{1},h_{0})).\]
\end{prop}
\vspace{.1in}
\underline{\em{Proof:}}

\vspace{.2in}
First notice that $((\lambda,\sigma),U)$ is a twisted action of $G\times H$ on
$A$: conditions a), b) and c) in  \cite[Def. 2.1]{pr} are easily checked, and,
for $(x_{0},h_{0})$, $(x_{1},h_{1})$, and $(x_{2},h_{2})$ in $G\times H$, we
have
\[(\lambda,\sigma)_{(x_{0},h_{0})}[U((x_{1},h_{1}), (x_{2},h_{2}))]U
((x_{0},h_{0}), (x_{1}x_{2},h_{1}h_{2}))=\]
\[=\lambda_{x_{0}}\sigma_{h_{0}}\lambda_{x_{1}}(u(x_{2},h_{1}))\lambda_{x_{0}}(u(x_{1}x_{2},h_{0}))=\lambda_{x_{0}x_{1}}(u(x_{2},h_{0}h_{1}))\lambda_{x_{0}}(u(x_{1},h_{0})=\]
\[=U((x_{0}x_{1},h_{0}h_{1}),(x_{2},h_{2}))U((x_{0},h_{0}),(x_{1},h_{1})).\]

We now construct maps
\[i_{A}:A\rightarrow {\cal M}(A\times_{\lambda}G\times_{\gamma^{\sigma
,u}}H)\mbox{\hspace{.1in}and \hspace{.1in}}i_{G\times H}:G\times
H\rightarrow{\cal UM}(A\times_{\lambda}G\times_{\gamma^{\sigma ,u}}H)\]
satisfying
\[i_{A}((\lambda,\sigma)_{(x,h)}(a))=i_{G\times H}(x,h)i_{A}(a)i_{G\times
H}(x,h)^{*}\mbox{\hspace{.1in}and \hspace{.1in}}\]
\[ i_{G\times H}(x_{0},h_{0})i_{G\times
H}(x_{1},h_{1}))=i_{A}(U((x_{0},h_{0}),(x_{1},h_{1}))i_{G\times
H}(x_{0}x_{1},h_{0}h_{1}),\]
for all $x_{i}\in G$, $h_{i}\in H$, and  $a\in A$.

If $\alpha$ is an action of a group $K$ on a C*-algebra B, $b\in {\cal M}(B)$,
and  $\mu$ is a bounded complex Radon measure with compact support on $G$, ,
let  $M(b,\mu)$ denote the multiplier of $B\times_{\alpha}K$ defined by
\[(M(b,\mu)f)(t)=b\int_{K}\alpha_{s}(f(s^{-1}t))d\mu (s),\]
fpr $f\in C_{c}(K,B)$.

Now define
\[i_{A}(a)=M(M(a,\delta_{1_{G}}),\delta_{1_{H}})\mbox{\hspace{.1in}and
\hspace{.1in}} i_{G\times H}(x,h)=M(M(1_{A},\delta_{x}),\delta_{h}),\]
where $\delta_{t}$ denotes the  point mass at $t$.

For $f\in C_{c}(G\times H, A)$, explicit formulas are given by:
\[(i_{A}(a)f)(x,h)=a f(x,h),\mbox{and}\]
\[(i_{G\times
H}(x_{0},h_{0})f)(x,h)=u^{*}(x_{0},h_{0})u(x,h_{0})\lambda_{x_{0}}\sigma_{h_{0}}(f(x_{0}^{-1}x,h_{0}^{-1}h)).\]
It follows that
\[(i_{G\times
H}^{*}(x_{0},h_{0})f)(x,h)=u(x,h_{0}^{-1})\sigma_{h_{0}^{-1}}\lambda_{x_{0}^{-1}}(f(x_{0}x,h_{0}h)).\]
The pair $(i_{A},i_{G\times H})$ is covariant:
\[(i_{G\times H}(x_{0},h_{0})i_{A}(a)i^{*}_{G\times H}(x_{0},h_{0})f)(x,h)=\]
\[=u^{*}(x_{0},h_{0})u(x,h_{0})\lambda_{x_{0}}\sigma_{h_{0}}[au(x^{-1}_{0}x,h_{0}^{-1})\sigma_{h_{0}^{-1}}\lambda_{x^{-1}_{0}}(f(x,h))]=\]
\[=(i_{A}(\lambda_{x_{0}}\sigma_{h_{0}}(a))f)(x,h),\]
and
\[(i_{G\times H}(x_{0},h_{0})i_{G\times H}(x_{1},h_{1}))(x,h)=\]
\[u^{*}(x_{0},h_{0})u(x,h_{0})\lambda_{x_{0}}\sigma_{h_{0}}[u^{*}(x_{1},h_{1})u(x_{0}^{-1}x,h_{1})\lambda_{x_{1}}\sigma_{h_{1}}(f(x_{1}^{-1}x_{0}^{-1}x,h_{1}^{-1}h_{0}^{-1}h))]=\]
\[=\lambda_{x_{0}}u(x_{1},h_{0})u^{*}(x_{0}x_{1},h_{0}h_{1})u(x,h_{0}h_{1})\lambda_{x_{0}x_{1}}\sigma_{h_{0}h_{1}}(f(x_{1}^{-1}x_{0}^{-1}x,h_{1}^{-1}h_{0}^{-1}h))=\]
\[=U((x_{0},h_{0}),(x_{1},h_{1}))i_{G\times
H}((x_{0}x_{1},h_{0}h_{1})f)(x,y).\]
We next show that for any covariant representation $(\Pi,V)$ of \mbox{$(A,
G\times H,(\lambda,\sigma), U)$} on a Hilbert space $\cal H$ there is an
integrated form $\Pi\times V$ on
\mbox{$A\times_{\lambda}G\times_{\gamma^{\sigma ,u}}H$}.
Let $V_{G}$ and $V_{H}$ be the restrictions of $V$ to $G$ and $H$,
respectively. Then $(\Pi,V_{G})$ is a covariant representation of
\mbox{$(A,G,\lambda)$} and, if $\Pi\times V_{G}$ denotes its integrated form,
then \mbox{$(\Pi\times V_{G},V_{H})$} is a covariant representation of
\mbox{$(A\times_{\lambda}G,H,\gamma^{\sigma ,u})$}.
So  $\Pi\times V_{G}\times V_{H}$ is a non-degenerate representation of
\mbox{$A\times_{\lambda}G\times_{\gamma^{\sigma ,u}}H$} and
\[\Pi=\Pi\times V_{G}\times V_{H}\circ i_{A}\mbox{ and }V=\Pi\times V_{G}\times
V_{H}\circ i_{G\times H}.\]
Finally, the set \mbox{$\{i_{A}\times i_{G\times H}(f):f\in L^{1}(G\times
H,A)\}$}, where
\[[i_{A}\times i_{G\times H}(f)](x,h)=\int_{G\times H}i_{A}[f(x,h)]i_{G\times
H}(x,h)d(x,y)\]
 is a dense subspace of \mbox{$A\times_{\lambda}G\times_{\gamma^{\sigma
,u}}H$}, which ends the proof.
\begin{flushright}
Q.E.D.
\end{flushright}
\section{The generalized fixed-point algebras.}
\label{meq}
  With the example of quantum Heisenberg manifolds in mind, we now discuss the
situation described in Section \ref{gen} in the case of some particular actions
$\lambda$ and $\sigma$ on a commutative C*-algebra $C_{0}(M)$. We prove that if
the action $\sigma$ is proper, then so is $\gamma^{\sigma,u}$  (in the sense of
\cite{rfpa}), and that if $\sigma$ is also free then $\gamma^{\sigma,u}$ is
saturated (\cite{rfpa}). Besides, for $\lambda$ and $\sigma$ free and proper,
the  generalized fixed-point algebras under $\gamma^{\sigma,u}$ and
$\gamma^{\lambda,u^{*}}$ respectively are strong-Morita equivalent.
\newline More specifically, we show that the space $C_{c}(M)$ can be made into
a dense submodule of an equivalence bimodule for the generalized fixed-point
algebras. Part of this is done by adapting to our situation the techniques of
{\cite[Situation 10]{rfsit}.
\vspace{.1in}
\newline
\underline{\bf{Assumptions and notation.}} Throughout this section $M$ denotes
a locally compact Hausdorff space, and $\beta M$ its Stone-Cech
compactification. The groups G and H act on $M$ by commuting actions $\lambda$
and $\sigma$, respectively.
In this context, if $T$ denotes the unit circle, the cocycle $u$ of Section
\ref{gen} consists of continuous functions \mbox{$u(x,h):M\rightarrow T$}, for
$(x,h)\in G\times H$, such that, for any $f\in C_{0}(M)$ the map
$(x,h)\rightarrow u(x,h)f$ is continuous for the supremum norm. As in Section
\ref{gen} we require the cocycle conditions:
\[u(x_{1}x_{2},h)=u(x_{1},h)\lambda_{x_{1}}u(x_{2},h)\mbox{\hspace{.2in}and\hspace{.2in}}u(x,h_{1}h_{2})=u(x,h_{1})\sigma_{h_{1}}u(x,h_{2}),\]
for $x, x_{i}\in G$ and $h,h_{i}\in H$. Notice that if these conditions are
satisfied for $u$ they also hold for $u^{*}$. We denote by $\gamma^{\sigma ,u}$
and $\gamma^{\lambda ,u^{*}}$ the actions of $H$ and $G$ on
$C_{0}(M)\times_{\lambda}G$ and $C_{0}(M)\times_{\sigma}H$ respectively, as
defined in Proposition \ref{act}.
\begin{prop}
\label{main}
In the notation above, if $\sigma$ is proper, so is the action $\gamma^{\sigma
,u}$ of $H$ on $C_{0}(M)\times_{\lambda}G$.
The generalized fixed-point algebra $D^{\sigma ,u}$ of
$C_{0}(M)\times_{\lambda}G$ under $\gamma^{\sigma ,u}$ consists of the closure
in ${\cal M}(C_{0}(M)\times_{\lambda}G)$ of the linear span of the set
$\{P_{\sigma ,u}(E^{*}*F):E,F\in C_{c}(M\times G)\}$, where $P_{\sigma ,u}$
denotes the linear map $P_{\sigma ,u}:C_{c}(M\times G)\rightarrow {\cal
M}(C_{0}(M)\times_{\lambda}G)$ defined by
\[(P_{\sigma ,u}(F))(m,x)=\int_{H}(\gamma_{h}^{\sigma ,u}(F))(m,x)dh,\]
for $F\in C_{c}(M\times G)$, and $(m,x)\in M\times G$.
\newline Furthermore, $P_{\sigma ,u}$ satisfies

i) $P_{\sigma ,u}(F^{*})=P_{\sigma ,u}(F)^{*}$

ii) $P_{\sigma ,u}(F)\geq 0$, for $F\geq 0$, where $F$ and $P_{\sigma ,u}(F)$
are viewed as elements of ${\cal M}(C_{0}(M)\times_{\lambda}G).$

iii) $P_{\sigma ,u}(F*\Phi)=P_{\sigma ,u}(F)*\Phi$ and $P_{\sigma
,u}(\Phi*F)=\Phi*P_{\sigma ,u}(F)$,
\newline for any $\Phi\in {\cal M}(C_{0}(M)\times_{\lambda} G)$ carrying
$C_{c}(M\times G)$ into itself and such that $\gamma^{\sigma
,u}_{h}(\Phi)=\Phi$ for any $h\in H$.

\end{prop}
\vspace{.1in}
\underline{\em {Proof:}}
\vspace{.2in}

We check conditions 1) and 2) of \cite[Def. 1.2]{rfpa}. Let $B=C_{c}(M\times
G)$. Then $B$ is a dense *-subalgebra of $C_{0}(M)\times_{\lambda}G$, and it is
invariant under $\gamma^{\sigma ,u}$.

We now show that, for $E,F\in B$, the map $h\rightarrow E*\gamma^{\sigma
,u}_{h}(F^{*})$ is in $L^{1}(H,C_{0}(M)\times_{\lambda}G)$. For $(m,x)\in
M\times G$ we have
\[[E*\gamma^{\sigma
,u}_{h}(F^{*})](m,x)=\int_{G}E(m,y)[u(y^{-1}x,h)](\lambda_{y^{-1}}m)\overline{F}(\lambda_{x^{-1}}\sigma_{h^{-1}}m,x^{-1}y)\Delta_{G}(x^{-1}y)dy.\]
Since $\sigma$ is proper and $supp(E)$ and $supp(F)$ are compact, then the set
\[\{h\in H : \sigma_{h^{-1}}\lambda_{x^{-1}}m\in supp_{M}(F) \mbox{ for
}(m,x)\in supp_{M}(E)\times supp_{G}(E)supp_{G}(F)^{-1}\}\]
is compact. Therefore $h\rightarrow E*\gamma^{\sigma ,u}_{h}(F^{*})$ and
$h\rightarrow \Delta_{H}^{-1/2}(h)E*\gamma^{\sigma ,u}_{h}(F^{*})$ are  in
$C_{c}(H,B)\subseteq L^{1}(H,{\cal M}(C_{0}(M)\times_{\lambda}G))$.
\newline For $F\in B$ and $m_{0}\in M$, let $N$ be a neighborhood of $m_{0}$
with compact closure. Then there exists a compact set $K$ in $H$ such that
\[ P_{\sigma ,u}(F)(m,x)=\int_{K}(\gamma^{\sigma ,u}_{h}F)(m,x)dh,\]
for all $(m,x)\in N\times G$, which shows that $P_{\sigma ,u}(F)$ is
continuous.
Since $supp_{G}(P_{\sigma ,u}(F))$ is compact, then $P_{\sigma ,u}(F)$ is
bounded on $supp_{M}(F)\times G$.
Besides, for all $(m,x)\in M\times G$ and $h\in H$, we have $|P_{\sigma
,u}F(m,x)|=|P_{\sigma ,u}F(\sigma_{h}m,x)|$, and
$supp_{M}(P_{\sigma,u}(F))\subset\sigma_{H}(supp_{M}(F)).$
\newline Therefore $P_{\sigma ,u}(F)\in C_{c}(\beta M\times G)\subseteq {\cal
M}(C_{0}(M)\times_{\lambda}G),$ and, as a multiplier, $P_{\sigma ,u}(F)$
carries $B$ into itself.
\newline Notice now that the fact that $h\rightarrow E*\gamma_{h}^{\sigma
,u}(F)$ is in $L^{1}(H,C_{0}(M)\times_{\lambda}G)$ implies that the integral
\mbox{$\int_{H}\gamma^{\sigma ,u}_{h}(F)dh$} makes sense as an integral in the
completion of ${\cal M}(C_{0}(M)\times_{\lambda}G)$, viewed as a locally convex
linear space, for the topology induced by the set of seminorms $\{\|\mbox{
}\|_{F}:F\in B\}$, where
\[\|\Phi\|_{F}=\|F*\Phi\|_{C_{0}(M)\times_{\lambda}G}+\|\Phi*F\|_{C_{0}(M)\times_{\lambda} G}\]
for $\Phi\in {\cal M}(C_{0}(M)\times_{\lambda}G).$

A straightforward application of Fubini's theorem shows that
\[\int_{H}(E*\gamma_{h}^{\sigma ,u}(F))(m,x)dh = (E*P_{\sigma ,u}(F))(m,x),\]
 for any $E$, $F\in B$, $(m,x)\in M\times G$, and it follows that
\[\int_{H}\gamma^{\sigma ,u}_{h}(F)dh=P_{\sigma ,u}(F),\]
in the sense mentioned above.
\newline Also, since the positive cone is closed, and involution and the
extension of $\gamma^{\sigma ,u}$ are continuous for the topology of ${\cal
M}(C_{0}(M)\times_{\lambda}G)$ defined above, $P_{\sigma, u}$ satisfies i),
ii), and iii) stated above.
\vspace{.1in}

Set now $<E,F>_{\sigma}=P_{\sigma ,u}(E^{*}*F)$, for $E,$ $F\in B$. We have
shown that $\gamma^{\sigma ,u}$ is proper. The generalized fixed-point algebra
$D^{\sigma ,u}$ (\cite[Def.1.4]{rfpa}) of $C_{0}(M)\times_{\lambda} G$ under
$\gamma^{\sigma ,u}$ consists of the closure in ${\cal
M}(C_{0}(M)\times_{\lambda}G)$ of the linear span of the set $\{
<E,F>_{\sigma}:E,F\in B\}.$
\begin{flushright}
Q.E.D.
\end{flushright}
\vspace{.2in}
\begin{lemma}
\label{approx}
Assume that $\sigma$ is proper and let $\{\Phi_{N,\epsilon,K}\}$ be a net in
\newline $C_{c}(M\times G\times H)$, indexed by decreasing neighborhoods $N$ of
$1_{G\times H}$, decreasing $\epsilon >0$, and increasing compact subsets $K$
of $M$, satisfying

i) $supp_{G\times H}(\Phi_{N,\epsilon,K})\subset N $

ii) $|\int_{G\times H}\Phi_{N,\epsilon, K}(m,x,h)dxdh-1|<\epsilon$, for all
$m\in K$

iii)  There exists a real number $R$ such that $\int_{G\times
H}|\Phi_{N,\epsilon,K}(m,x,h)|dxdh\leq R$, for all $m\in K$, and for all $K$,
$\epsilon$ and $N$.

Then $\{\Phi_{N,\epsilon,K}\}$ is an approximate identity for $C_{c}(M\times
G\times H)\subset \newline C_{0}(M)\times_{\lambda}G\times_{\gamma^{\sigma
,u}}H$ in the inductive limit topology.
\end{lemma}
\vspace{.1in}
\underline{Proof:}

\vspace{.1in}
Let $\psi\in C_{c}(M\times G\times H)$ and $\delta >0$ be given. Then
\[|(\Phi_{N,\epsilon,K}*\Psi-\Psi)(m,x,h)|\leq\]
\[\leq |\int_{H\times
G}[u^{*}(y,k)(m)u(x,k)(m)-1]\Phi_{N,\epsilon,K}(m,y,k)\Psi(\sigma_{k^{-1}}\lambda_{y^{-1}}m,y^{-1}x, k^{-1}h)dkdy|+\]
\[+|\int_{H\times G}\Phi_{N,\epsilon,K}(m,y,k)dydk-1||\Psi(m,x,h)|+\]
\[+|\int_{H\times
G}\Phi_{N,\epsilon,K}(m,y,k)[\Psi(\sigma_{k^{-1}}\lambda_{y^{-1}}m,y^{-1}x,k^{-1}h)-\Psi(m,x,h)]dydk|\leq\delta,\]
for appropriate choices of $\epsilon$ and $N$.
\begin{flushright}
Q.E.D.
\end{flushright}

\begin{prop}
\label{sat}
 If the action $\sigma$ is free and proper, then  $\gamma^{\sigma ,u}$ is
saturated.
\end{prop}
\vspace{.1in}
\underline{Proof:}

\vspace{.1in}
Let $J$ denote the ideal of $C^{*}_{r}(H, C_{0}(M\times_{\lambda}G))$
consisting of maps $h\mapsto
\Delta_{H}^{-1/2}(h)E*\gamma_{h}^{\sigma,u}(F^{*})$, for $E,F\in C_{c}(M\times
G)$.
In order to show that $J$ is dense in $C^{*}_{r}(H,
C_{0}(M)\times_{\lambda}G))$ we prove that $J$ contains an approximate identity
for $C_{c}(M\times G\times H)$.

Let $N,\epsilon$, and $K$ as in Lemma \ref{approx} be given. We assume without
loss of generality that the closure of $N$ is compact. Fix an open set $U$ with
compact closure such that $K\subset U$. Choose neighborhoods $N_{G}$ and
$N_{H}$ of $1_{G}$ and $1_{H}$, respectively, such that $N_{G}\times
N_{H}\subset N$, $|\Delta_{G}(x)-1|<\epsilon_{1}$ for all $x\in N_{G}$ and
$|u^{*}(y,h)(m)u(x,h)(m)-1|<\epsilon_{2}$, for all  $h\in N_{H},m\in U$, $x,
y\in V$, $V$ being a fixed open set with compact closure containing $N_{G}$,
and for some $\epsilon_{1}$ and $\epsilon_{2}$ to be chosen later.

The action of $G\times H$ on $M\times G$ defined by
$(x,h)(m,y)=(\lambda_{x}\sigma_{h}m,xy)$ is free and proper, so for each
$(m,y)\in K\times \overline{N_{G}}$ we can choose (\cite[Situation 10]{rfsit})
a neighborhood $U_{(m,y)}\subset U\times V$ of $(m,y)$ such that
\[\{(x,h): (x,h)(U_{(m,y)})\cap U_{(m,y)}\not =\emptyset \}\subset N_{G}\times
N_{H}.\]
 Take a finite subcover $\{U_{1}, U_{2},...,U_{n}\}$ of
$\{U_{(m,y)}\}_{(m,y)\in K\times \overline{N_{G}}}$ and, for each $i=1,...,n$,
let $F_{i}\in C^{+}_{c}(M\times G)$ be such that $supp(F_{i})\subset U_{i}$,
and $\int_{G}\sum_{i}F_{i}(m,x)dx=1$ for all $m\in K$.

Now we can find (\cite[Situation 10]{rfsit}) functions $G_{i}\in
C_{c}^{+}(M\times G)$ such that $supp(G_{i})\subset supp(F_{i})$, and
\[|F_{i}(m,y)-G_{i}(m,y)\int_{G\times
H}G_{i}(\lambda_{x^{-1}}\sigma_{h^{-1}}m,x^{-1}y)dxdh|<\epsilon_{3},\]
for all $(m,y)\in M\times G$, and some $\epsilon_{3}$ to be chosen later.

Now set
\[\Phi_{N,\epsilon,K}(m,x,h)=\sum_{i}\Delta_{H}^{-1/2}(h)[G_{i}*\gamma_{h}^{\sigma,u}(G_{i}^{*})(m,x).\]
Then, \[|\int_{H\times G}\Phi_{N,\epsilon,K}(m,x,h)dxdh-1|=\]
\[=\sum_{i}\int_{G\times G\times
H}\Delta_{G}(x^{-1}y)[u^{*}(y,h)u(x,h)](m)G_{i}(m,y)G(\lambda_{x^{-1}}\sigma_{h^{-1}}m,x^{-1}y)dxdydh -\]
\[-\sum_{i}\int_{G}F_{i}(m,y)dy|\leq\]
\[\leq
|\sum_{i}\int_{V}([u^{*}(y,h)(m)u(x,h)(m)\Delta_{G}(x^{-1}y)-1]G_{i}(m,y)\int_{G\times H}G_{i}(\lambda_{x^{-1}}\sigma_{h^{-1}}m,x^{-1}y)dxdh)dy|+\]
\[+|\sum_{i}\int_{V}G_{i}(m,y)\int_{G\times
H}G_{i}((\lambda_{x^{-1}}\sigma_{h^{-1}}m,x^{-1}y)dxdh-F_{i}(m,y)dy|<\epsilon,\]
for appropriate choices of $\epsilon_{1}$,  $\epsilon_{2}$, and $\epsilon_{3}$.
\newline Besides, $supp(\Phi_{N,\epsilon,K})\subset N_{G}\times N_{H}\subset
N$.
Finally, a similar argument shows that from some $N_{0}$ and $\epsilon_{0}$ on
we have
\[\int_{H\times G}|\Phi_{N, \epsilon, K}(m,x,h)|dxdh\leq R,\]
for some real number $R$, and all $m\in K$.

\begin{flushright}
Q.E.D.
\end{flushright}
\vspace{.1in}
\underline{{\bf Assumptions.}} We next compare the generalized fixed-point
algebras obtained when the roles of $\sigma$ and $\lambda$ are reversed. That
is why we require symmetric conditions on these two actions. So, we assume from
now on that both $\lambda$ and $\sigma$ are free and proper actions.
\vspace{.1in}
\newline \underline{\bf Notation.}  Let $C^{\sigma,u}$ denote the subalgebra of
${\cal M}(C_{0}(M)\times_{\lambda}G)$ consisting of functions $\Phi\in
C_{c}(\beta M\times G)$ such that the projection of $supp_{M}(\Phi)$ on $M/H$
is precompact and $\gamma_{h}^{\sigma,u}\Phi=\Phi$ for all $h\in H$.
\vspace{.1in}
\begin{rk}
\label{csigma}
Notice that, for $F\in C_{c}(M\times G)$, we have that
$supp_{M}(P_{\sigma,u}F)\subset \sigma_{H}(supp_{M}(F))$, and therefore
$C^{\sigma,u}$ contains the image of $P_{\sigma,u}$.
\end{rk}

\vspace{.1in}
\begin{lemma}
\label{aid}
Let $\{\Phi_{N,\epsilon}\}$ be a net in $C^{\sigma,u}$, indexed by decreasing
neighborhoods $N$ of $1_{G}$, increasing compact subsets $K$ of $M$, and
decreasing $\epsilon >0$, and such that

1) $supp_{G}(\Phi_{N,\epsilon,K})\subseteq N$

2) $|\int_{G}\Delta_{G}^{1/2}(x) \Phi_{N,\epsilon}(m,x)dx-1|<\epsilon$ for all
$m\in K$.

3) There is a real number $R$ such that $\int_{G}|\Phi_{N,\epsilon}(m,x)|dx\leq
R$, for all $m\in K$, and for all $N$ and $\epsilon$ from some $N_{0}$ and
$\epsilon_{0}$ on.

Then $\{\Phi_{N,\epsilon,K}\}$ is an approximate identity for $C^{\sigma ,u}$.
\end{lemma}
\vspace{.1in}
\underline{\em {Proof:}}
\vspace{.1in}

Let $\Psi\in C^{\sigma,u}$ and $\delta>0$ be given. Fix a neighborhood $N'$ of
$1_{G}$ with compact closure, and let $K'\subset M$ be a compact set such that
$\Pi_{H}(supp_{M}\Psi)\subset\Pi_{H}(K')$, where $\Pi_{H}$ denotes the
canonical projection on $M/H$.

As in Lemma \ref{approx}, we can find $N_{0}\subset N'$, $\epsilon_{0}$, and
$K_{0}$ such that, from $N_{0}$, $\epsilon_{0}$, and $K_{0}$ on, we have
\[|(\Phi_{N,\epsilon,K}*\Psi-\Psi)(m,x)|<\delta,\]
for all $m\in \lambda_{N'}(K')$.

Therefore, if $m\in supp(\Phi_{N,\epsilon,K}*\Psi-\Psi)$, then we have that
$\sigma_{h}m\in \lambda_{N'}(K')$, for some $h\in H$.
On the other hand we have that
\[|(\Phi_{N,\epsilon,K}*\Psi-\Psi)(\sigma_{h}m,x)|=|(\Phi_{N,\epsilon,K}*\Psi-\Psi)(m,x)|,\]
for all $h\in H$, $m\in M$, and $x\in G$, because $\Phi_{N,\epsilon,K}$ and
$\Psi\in C^{\sigma,u}$.
  This shows that $|(\Phi_{N,\epsilon,K}*\Psi-\Psi)(m,x)|<\delta$ for all $m\in
M$. Therefore $\Phi_{N,\epsilon,K}*\Psi$ converge to $\Psi$ in the multiplier
algebra norm.

\begin{flushright}
Q.E.D.
\end{flushright}

\begin{rk}
\label{raid}
Notice that  Lemma \ref{aid} above also holds, with a similar proof, if
condition 2) is  replaced by

2') $|\int_{G} \Phi_{N,\epsilon,K}(m,x)dx-1|<\epsilon$ for all $m\in K$.
\end{rk}
\vspace{.2in}
\begin{prop}
\label{imp}
The generalized fixed point algebra $D^{\sigma ,u}$ is the closure in ${\cal
M}(C_{0}(M)\times_{\lambda}G)$ of $C^{\sigma ,u}$.
\end{prop}
\vspace{.1in}
\underline{{\em Proof:}}
\vspace{.1in}

In view of property iii) in Proposition \ref{main}, it suffices to show that
the span of the set
\newline \centerline {$\{P_{\sigma ,u}(E^{*}*F):E,F\in C_{c}(M\times G)\}$}
\newline contains an approximate identity for $C^{\sigma ,u}$.
\newline For a given compact set $K\subset M$, let us fix an open set $U$ of
compact closure containing $K$. Then the set \mbox{$L=\{h\in H: \sigma_{h}m\in
\overline{U}\mbox{ for some }m\in K\}$} is compact.

Let N be a given neighborhood of $1_{G}$ and $\epsilon >0$.
As in \cite[Sit. 10, first lemma]{rfsit}, we can take an open cover
$\{U_{1},U_{2},...,U_{n}\}$ of $K$, such that $U_{i}\subseteq U$ and $U_{i}\cap
\lambda_{x}U_{i}\not = \emptyset$ only if $x\in N$.
For each $i=1,...,n$, let $H_{i}\in C_{c}^{+}(M\times G)$ be such that
$supp(H_{i})\subset U_{i}\times N$, and $\sum_{i}H_{i}$ is strictly positive on
$K\times 1_{G}$.
Then $\sum_{i}\int_{H\times G} H_{i}(\sigma_{h^{-1}}m,y)dhdy>0$ for all
$m\in K$. Therefore, we can find functions $F_{i}\in C_{c}^{+}(M\times G)$
such that $supp(F_{i})\subset supp(H_{i})$ and \mbox{$\int_{H\times
G}F_{i}(\sigma_{h^{-1}}m,y)dhdy$}$=1$ for all $m\in K$.
Now, the action of $G$ on $M\times G$ given by
$\alpha_{x}(m,y)=(\lambda_{x}m,xy)$ is free and proper, so the second lemma in
\cite[Situation 10]{rfsit} applies and for each $i=1,...,n$ we can find
$G_{i}\in C^{+}_{c}(M\times G)$ such that $supp(G_{i})\subseteq supp(F_{i})$
and
\[|F_{i}(m,y)-G_{i}(m,y)\int_{G}G_{i}(\lambda_{x^{-1}}m,x^{-1}y)dx|<\delta/n,\]
for all $m\in M, y\in G$, and some positive number $\delta$ to be chosen later.
Set now $\Phi_{N,\epsilon,K}=\sum_{i=1}^{i=n}P_{\sigma ,u}(G_{i}*J_{i})$, where
$J_{i}(m,x)=G_{i}(\lambda_{x^{-1}}m,x^{-1})$.
We have
\[\Phi_{N,\epsilon,K}(m,x)=\sum_{i}\int_{H}u(x,h)\int_{G}G_{i}(\sigma_{h^{-1}}m,y)G_{i}(\sigma_{h^{-1}}\lambda_{x^{-1}}m,x^{-1}y)dy,\]
so, since $supp(G_{i})\subseteq supp(F_{i})$, it follows that
$supp_{G}(\Phi_{N,\epsilon,K})\subseteq N$.
\newline Besides, if $m\in K$,
\[|\int_{G}\Phi_{N,\epsilon,K}(m,x)dx-1|=\]
\[=
|\sum_{i}\int_{H}\int_{G}[u(x,h)G_{i}(\sigma_{h^{-1}}m,y)\int_{G}G_{i}(\sigma_{h^{-1}}\lambda_{x^{-1}}m,x^{-1}y)dx-F_{i}(\sigma_{h^{-1}}m,y])dy dh |<\epsilon,\]
for a suitable choice of $\delta$, if $N$ is chosen to have $|u(x,h)-1|$ small
enough for all $x\in N$ and $ h\in L$.
\newline Finally, from some $\epsilon_{0}$ and $N_{0}$ on,
$\int_{G}|\Phi_{N,\epsilon, K}(m,x)|dx\leq R$, for some real number $R$ and all
$m\in K$.
\newline Then, by Remark \ref{raid}, $\{\Phi_{N,\epsilon,K}\}$ is an
approximate identity for $C^{\sigma ,u}$.
\begin{flushright}
Q.E.D.
\end{flushright}
\vspace{.1in}

We will later make use of the following variation of the construction in the
proof of Theorem \ref{imp}.
\vspace{.1in}
\begin{rk}
\label{fg}
The span of the set
\[\{P_{\sigma
,u}(F):F(m,x)=\Delta^{-1/2}(x)e_{i}(m)\overline{e}_{i}(\lambda_{x^{-1}}m) ,
e\in C_{c}(M)\}\]
contains an approximate identity for $C{\sigma ,u}$.
\end{rk}
\vspace{.1in}
\underline{\em Proof:}

\vspace{.1in}
In the notation of Proposition \ref{imp}, let $\{f_{i}\}\subset C_{c}^{+}(M)$
be such that $supp(f_{i}\subset U_{i}$, and
$\int_{H}\sum_{i}f_{i}(\sigma_{h^{-1}}m)>0$, for all $m\in K$. Since the action
$\lambda$ is proper we can get $g_{i}\in C^{+}_{c}(M)$ such that
$supp(g_{i})\subseteq supp(f_{i})$ and
\newline $|f_{i}(m)-g_{i}(m)\int_{G}g_{i}(\lambda_{x^{-1}}m)dx|<\delta$ for all
$m\in M$ and a given positive number $ \delta.$
Then, if we let
$L_{i}(m,x)=\Delta_{G}^{-1/2}(x)g_{i}(m){g_{i}}(\lambda_{x^{-1}}m)$ we have
that, for an appropriate choice of $\delta$ in terms of $\epsilon$, the
function $\Phi_{N,\epsilon,K}=\sum_{i}P_{\sigma,u}(L_{i})$ can be shown (by an
argument quite similar to that in Proposition \ref{imp}) to satisfy the
hypotheses of Lemma \ref{aid}.
\begin{flushright}
Q.E.D.
\end{flushright}
\vspace{.2in}
\underline{\bf Notation.}
We denote by $_{\lambda}<$ , $>$ and $<$ , $>_{\lambda}$ the $C_{c}(M\times
G)$-valued maps defined on $C_{c}(M)\times C_{c}(M)$ by
\[_{\lambda}<f,g>(m,x)=\Delta^{-1/2}_{G}(x)f(m)\overline{g}(\lambda_{x^{-1}}m)\]\[\mbox{\hspace{.1in}and\hspace{.1in}}<f,g>_{\lambda}(m,x)=\Delta^{-1/2}_{G}(x)\overline{f}(m)g(\lambda_{x^{-1}}m),\]
where $f,g\in C_{c}(M).$
\vspace{.2in}
\begin{rk}
\label{bc}
It is a well known result (\cite[Situation 2]{rfsit}) that $C_{c}(M)$ is a left
(resp. right) $C_{c}(M\times G)$-rigged module for $_{\lambda}<$ , $>$ (resp.
$<$ , $>_{\lambda}$) and the actions given by:
\[(\Phi\cdot f)(m)=\int_{G}\Delta^{1/2}_{G}(y)\Phi(m,y)f(\lambda_{y^{-1}}m)dy\]
\[\mbox{\hspace{.1in}and\hspace{.1in}}(f\cdot\Phi)(m)=\int_{G}\Delta_{G}^{-1/2}(y)\Phi(\lambda_{y^{-1}}m,y^{-1})f(\lambda_{y^{-1}}m)dy,\]
for $\Phi\in C_{c}(M\times G)$
It is easily checked that, by taking  $\Phi\in C_{c}(\beta M\times G)$ in the
formulas above, one makes $C_{c}(M)$ into a $C_{c}(\beta M\times G)$-module
with inner product. Of course it is no longer a rigged space because the
condition of density fails.
\end{rk}
\vspace{.2in}
\begin{prop}
\label{rigged}
Let $C^{\sigma ,u}\subseteq C_{c}(\beta M\times G)$ act on $C_{c}(M)$ on the
left and on the right as in Remark \ref{bc}. For $f,g\in C_{c}(M)$ define
\[<f,g>_{D^{\sigma ,u}}=P_{\sigma
,u}(<f,g>_{\lambda})\mbox{\hspace{.1in}and\hspace{.1in}}_{D^{\sigma
,u}}<f,g>=P_{\sigma ,u}(_{\lambda}<f,g>).\]
Then $C_{c}(M)$ is a left (resp. right) $C^{\sigma ,u}$-rigged space with
respect to \mbox{$_{D^{\sigma ,u}}<$ , $>$} (resp. $<$ , $>_{D^{\sigma ,u}})$.
\end{prop}
\vspace{.1in}
\underline{{\em Proof:}}

\vspace{.1in}
The density condition follows from Remark \ref{fg}. All other properties follow
immediately from the fact that $_{\lambda}<$ , $>$ and $<$ , $>_{\lambda}$ are
inner products and from Remark \ref{csigma} and properties i), ii), and iii) of
$P_{\sigma ,u}$ shown in Proposition \ref{main}.

\begin{flushright}
Q.E.D.
\end{flushright}
We are now ready to show the main result of this section.

\vspace{.2in}
\begin{thm}
\label{dmeq}
Let $\lambda$ and $\sigma$ be free and proper commuting actions of locally
compact groups $G$ and $H$ respectively on a locally compact space $M$. Let $u$
be a cocycle as in Proposition \ref{act}. Then the generalized fixed-point
algebras $D^{\sigma ,u}$ and $D^{\lambda ,u^{*}}$ of the actions
$\gamma^{\sigma ,u}$ and $\gamma^{\lambda ,u^{*}}$ on
$C_{0}(M)\times_{\lambda}G$ and \mbox{$C_{0}(M)\times_{\sigma}H$,}
respectively, are strong-Morita equivalent.
\end{thm}
\vspace{.1in}
\underline{\em Proof:}

\vspace{.1in}
By Proposition \ref{rigged}, $C_{c}(M)$ is a left $C^{\sigma ,u}$-rigged space
and a right $C^{\lambda ,u^{*}}$-rigged space under
\[ (\Phi\cdot
f)(m)=\int_{G}\Delta_{G}^{1/2}(y)\Phi(m,y)f(\lambda_{y^{-1}}m)dy\mbox{\hspace{.1in},\hspace{.1in}}_{D^{\sigma ,u}}<f,g>=P_{\sigma ,u}(_{\lambda}<f,g>),\]
\[(f\cdot\Psi)(m)=\int_{H}\Delta_{H}^{-1/2}(h)\Psi(\sigma_{h^{-1}}m,h^{-1})f(\sigma_{h^{-1}}m)dh,\]
\[\mbox{\hspace{.1in}and\hspace{.1in}}<f,g>_{D^{\lambda ,u^{*}}}=P_{\lambda
,u^{*}}(<f,g>_{\sigma}),\]
where $f,g\in C_{c}(M)$, $\Phi\in C^{\sigma ,u}$ and $\Psi\in C^{\lambda
,u^{*}}$.
\newline Then $C_{c}(M)$ is an $C^{\sigma ,u}$-$C^{\lambda ,u^{*}}$ bimodule:
for $\Phi,$ $\Psi$ and $f$ as above we have
\[[(\Phi\cdot f)\cdot \Psi](m)=\]
\[\int_{H}\int_{G}\Delta_{H}^{-1/2}(h)\Delta_{G}^{1/2}(y)\Psi(\sigma_{h^{-1}}m,h^{-1})\Phi(\sigma_{h^{-1}}m,y)f(\sigma_{h^{-1}}\lambda_{y^{-1}}m)dydh=\]
\[=\int_{H}\int_{
G}\Delta_{H}^{-1/2}(h)\Delta_{G}^{1/2}(y)\Psi(\sigma_{h^{-1}}\lambda_{y^{-1}}m,h^{-1})\Phi(m,y)f(\sigma_{h^{-1}}\lambda_{y^{-1}}m)dydh=\]
\[=[\Phi\cdot(f\cdot\Psi)](m).\]
Besides, for $e,f,g\in C_{c}(M)$, we have
\[(_{D^{\sigma ,u}}<e,f>\cdot
g)(m)=\int_{G}\int_{H}u(y,h)e(\sigma_{h^{-1}}m)\overline{f}(\lambda_{y^{-1}}\sigma_{h^{-1}}m)g(\lambda_{y^{-1}}m)dhdy=\]
\[=(e<f,g>_{D^{\lambda ,u^{*}}})(m).\]
We now prove the continuity of the module structures with respect to the inner
products.
\newline Fix a measure $\mu$ of full support on $M$. Then, by \cite[6.1]{ph}
and \cite[7.7.5]{pd}, we have faithful representations $\Pi$ of $C^{\sigma ,u}$
on $L^{2}(M\times G)$ and $\Theta$ of $Im(P_{\lambda ,u^{*}}$ on $L^{2}(M\times
H)$ given by
\[(\Pi_{\Phi}\xi)(m,x)=\int_{G}\Phi(\lambda_{x}m,y)\xi(m,y^{-1}x)dx,\]
\[\mbox{
and}(\Theta_{\Psi}\eta)(m,h)=\int_{H}\Psi(\sigma_{h}m,k)\eta(m,k^{-1}h)dk,\]
where $\Phi\in C^{\sigma ,u}$, $\Psi\in C^{\lambda ,u^{*}}$, $\xi\in
L^{2}(M\times G)$ and $\eta\in L^{2}(M\times H)$.

Now, for $f\in C_{c}(M)$ and $\eta\in L^{2}(M\times H)$
\[<\Theta_{<f,f>_{D^{\lambda ,u^{*}}}}\eta,\eta>_{L^{2}(M\times H)}=\]
\[\int_{M\times G\times H\times
H}\sigma_{h^{-1}}(u^{*}(y,k))\Delta_{H}^{-1/2}(k)\overline{f}(\lambda_{y^{-1}}\sigma_{h}m).\]
\[.{f}(\lambda_{y^{-1}}\sigma_{k^{-1}h})\eta(m,k^{-1}h)\overline{\eta}(m,h)dkdhdydm=\]
\[=\|\xi(f,\eta)\|^{2}_{L^{2}(M\times G)},\]
where $\xi(f,\eta)\in L^{2}(M\times G)$ is given by
\[(\xi(f,\eta))(m,x)=\int_{H}u^{*}(x,h^{-1})\Delta_{H}^{-1/2}(h)f(\lambda_{x^{-1}}\sigma_{h}m)\eta(m,h)dh.\]
Then, if $\Phi\in C^{\sigma ,u}$
\[[\xi(\Phi\cdot f,\eta)](m,x)=\]
\[=\int_{G}\int_{H}u^{*}(x,h^{-1})\Delta_{H}^{-1/2}(h)\Delta_{G}^{1/2}(y)\Phi(\lambda_{x^{-1}}\sigma_{h}m,y)f(\lambda_{y^{-1}x^{-1}}\sigma_{h}m)\eta(m,h)dhdy=\]
\[=(U\Pi_{\Phi}U\xi(f,\eta))(m,x),\]
where $U$ denotes the unitary operator on $L^{2}(M\times G)$ defined by
$(U\xi)(m,x)=\Delta_{G}^{-1/2}(x)\xi(m,x^{-1})$.
\newline Thus we have
\[<\Theta_{<\Phi\cdot f,\Phi\cdot f>_{D^{\lambda
,u^{*}}}}\eta,\eta>_{L^{2}(M\times H)}=\|\xi(\Phi\cdot
f,\eta)|^{2}=\|U\Pi_{\Phi}U\xi(f,\eta)\|^{2}\leq\]
\[\leq \|\Phi\|^{2}\|\xi(f,\eta)\|^{2}=\|\Phi\|^{2}<\Theta_{<f,f>_{D^{\lambda
,u^{*}}}}\eta,\eta>_{L^{2}(M\times H)},\]
and it follows that
\[<\Phi\cdot f,\Phi\cdot f>_{D^{\lambda ,u^{*}}}\leq
\|\Phi\|^{2}<f,f>_{D^{\lambda ,u^{*}}},\]
as elements of $D^{\lambda ,u^{*}}$.
Analogously, one shows that, for $f\in C_{c}(M)$ and  $\xi\in L^{2}(M\times G)$
\[<\Pi_{_{D^{\sigma ,u}}<f,f>}\xi,\xi>_{L^{2}(M\times G)}|\eta(f,\xi)\|^{2},\]
for some $\eta(f,\xi)\in L^{2}(M\times H)$, and that, for $\Psi\in C^{\lambda
,u^{*}}$ one has
\[\eta(f\cdot\Psi,\xi)=(V\Theta_{\Psi^{*}}V)(\eta(f,\xi)),\]
where $V$ denotes the unitary operator in $L^{2}(M\times H)$ defined by
\newline $(V\eta)(m,h)=\Delta_{H}^{-1/2}(h)\eta(m,h^{-1})$. It follows that
\[_{D^{\sigma ,u}}<f\cdot \Psi,f\cdot\Psi>\leq \|\Psi\|^{2}_{D^{\sigma
,u}}<f,f>,\]
as elements of $D^{\sigma ,u}$.

Thus, we have proven that $C_{c}(M)$ is an $C^{\sigma ,u}-C^{\lambda ,u^{*}}$
equivalence bimodule. Now, if we define on $C_{c}(M)$ the norms
\[\|f\|^{2}_{D^{\sigma ,u}}=\|_{D^{\sigma
,u}}<f,f>\|\mbox{\hspace{.1in}and\hspace{.1in}}\|f\|^{2}_{D^{\lambda
,u^{*}}}=\|<f,f>_{D^{\lambda ,u^{*}}}\|,\]
it follows from \cite[3.1]{rfbl} that $\|\mbox{ }\|_{D^{\sigma ,u}}=\|\mbox{
}\|_{D^{\lambda ,u^{*}}}$ and that the completion of $C_{c}(M)$ with respect to
this norm gives, by continuity, an equivalence bimodule between $D^{\sigma ,u}$
and $D^{\lambda ,u^{*}}$.
\begin{flushright}
Q.E.D.
\end{flushright}
\vspace{.1in}
\begin{clly}
With the assumptions of Theorem \ref{dmeq}, the algebras
\newline $C^{*}_{r}(H, C_{0}(M)\times_{\lambda}G)$ and $C^{*}_{r}(G,
C_{0}(M)\times_{\sigma}H)$ are strong-Morita equivalent.
\end{clly}
\vspace{.1in}
\underline{\em Proof:}
\vspace{.1in}

The proof follows from Proposition \ref{sat}, Theorem \ref{dmeq}, and
\cite[1.7]{rfpa}.
\begin{flushright}
Q.E.D.
\end{flushright}

\section{Applications to quantum Heisenberg manifolds.}
In this section we apply the previous results to the computation of the
K-groups of the quantum Heisenberg manifolds. We recall the basic results and
definitions concerning those algebras. We refer the reader to \cite{rfhm} for
further details.

For each positive integer $c$, the Heisenberg manifold $M_{c}$ consists of the
quotient $G/D_{c}$, where $G$ is the Heisenberg group
\[G=\left\{\left( \begin{array}{ccc}
             1 & y & z \\
             0 & 1 & x \\
             0 & 0 & 1
             \end{array}
             \right) ; \mbox{ for real numbers }x,y,z \right\}\]
and $D_{c}$ is the discrete subgroup obtained when $x$, $y$ and $cz$ above are
integers.

The set of non-zero Poisson brackets on $M_{c}$ that are invariant under the
action of $G$ by left translation can be parametrized by two real numbers $\mu$
and $\nu$, with $\mu^{2}+\nu^{2}\not = 0$. A deformation quantization
$\{D^{c,\hbar}_{\mu\nu}\}_{\hbar\in R}$ of $M_{c}$ in the direction of a given
invariant Poisson bracket $\Lambda_{\mu\nu}$ was constructed in \cite{rfhm}.

The algebra $D^{c,\hbar}_{\mu\nu}$ can be described as a generalized
fixed-point algebra as follows. Let $M=R\times T$ and $\lambda^{\hbar}$ and
$\sigma$ be the commuting actions of $Z$ on $M$ induced by the homeomorphisms
\[\lambda^{\hbar}(x,y)=(x+2\hbar\mu,
y+2\hbar\nu)\mbox{\hspace{.1in}and\hspace{.1in}}\sigma(x,y)=(x-1,y).\]

Consider the action $\rho$ of $Z$ on $C_{0}(R\times
T)\times_{\lambda^{\hbar}}Z$ given by
\[(\rho_{k}\Phi)(x,y,p)=e(ckp(y-\hbar p\nu))\Phi(x+k,y,p),\]
where $e(x)=exp(2\pi ix)$ for any real number $x$.
The action $\rho$ defined above corresponds to the action $\rho$ defined in
\cite[p.539]{rfhm}, after taking Fourier transform in the third variable to get
the algebra denoted in that paper by $A_{\hbar}$ , and viewing $A_{\hbar}$ as a
dense *-subalgebra of $C_{0}(R\times T)\times_{\lambda^{\hbar}}Z$ via the
embedding $J$ defined in \cite[p.547]{rfhm}.

Notice that, for $M=R\times T$, $G=H=Z$, and $\hbar\not =0$, the actions
$\lambda^{\hbar}$ and $\sigma$ satisfy the hypotheses of Section \ref{meq} and
that the action $\rho$ defined above corresponds, in that context, to the
action we denoted by $\gamma^{\sigma ,u}$, where $u:Z\times Z\rightarrow {\cal
ZUM}(C_{0}(R\times T)$ is the 2-cocycle defined by
\[u(p,k)=e(ckp(y-\hbar p\nu)),\]
for $p,k\in Z$. Besides, \cite[Theorem 5.4]{rfhm} shows that the algebra
$D^{c,\hbar}_{\mu\nu}$ is the generalized fixed-point algebra of $C_{0}(R\times
T)\times_{\lambda^{\hbar}}Z$ under the action $\rho$, and it follows from the
proof of that theorem that $D^{c,\hbar}_{\mu\nu}$ is the algebra that we
denote, in the context of Section \ref{meq}, by $D^{\sigma ,u}$.

\begin{rk}
\label{dzero}

We will also use the fact that the algebra $\tilde{D}^{c,\hbar}_{\mu\nu}$
consisting of functions $\Phi\in C_{c}(\beta(R\times T)\times Z)$ satisfying
$\rho_{k}(\Phi)=\Phi$ for all $k\in Z$ is a dense *-subalgebra of
$D_{\mu\nu}^{c,\hbar}$. This follows from Remark \ref{csigma}, Proposition
\ref{imp}, and from the fact that $R\times T/\sigma$ is compact.

\end{rk}
\vspace{.1in}
\begin{thm}
\label{kth}
For $\hbar\not =0$ the K-groups of $D^{c,\hbar}_{\mu\nu}$ do not depend on
$\hbar$.

\end{thm}
\vspace{.1in}
\underline{\em Proof:}

\vspace{.1in}
It follows from Theorem \ref{dmeq} that, for $\hbar\not =0$,
$D^{c,\hbar}_{\mu\nu}$ is strong-Morita equivalent to the generalized
fixed-point algebra $E^{c,\hbar}_{\mu\nu}$ of $C_{0}(R\times
T)\times_{\sigma}Z$ under the action $\gamma^{\lambda^{\hbar}}$ of $Z$ defined
by
\[(\gamma^{\lambda^{\hbar}}_{p}\Phi)(x,y,k)=e(-ckp(y-\hbar p\nu))\Phi(x-2p\hbar
\mu,y-2p\hbar\nu,k).\]
Now, by Proposition \ref{sat}, $\gamma^{\lambda^{\hbar}}$ is saturated, so we
have (\cite[Corollary 1.7]{rfpa}) that $D^{c,\hbar}_{\mu\nu}$ is strong-Morita
equivalent to $C_{0}(R\times
T)\times_{\sigma}Z\times_{\gamma^{\lambda^{\hbar}}}Z$.

Besides, $\hbar\mapsto \lambda^{\hbar}$ is a homotopy between the
$\lambda^{\hbar}$'s, which shows (\cite[10.5.2]{bl}) that the K-groups of
$C_{0}(R\times T)\times_{\sigma}Z\times_{\gamma^{\lambda^{\hbar}}}Z$ do not
depend on $\hbar$.
On the other hand, since strong-Morita equivalent separable C*-algebras are
stably isomorphic (\cite{bgr}) and therefore have the same K-groups, we have
proven that the K-groups of $D^{c,\hbar}_{\mu\nu}$, for $\hbar\not =0$, do not
depend on $\hbar$.
\begin{flushright}
Q.E.D.
\end{flushright}
\vspace{.1in}
\underline{\bf Notation.} Since the algebras $D^{c,\hbar}_{\mu\nu}$ and
$D^{c,1}_{\hbar\mu,\hbar\nu}$ are isomorphic, we drop from now on  the constant
$\hbar$ from our notation  and absorb it into the parameters $\mu,\nu$.

\begin{rk}
\label{plusk}
Notice that, since for any pair of integers $k$ and $l$ the algebras
$D_{\mu\nu}^{c}$ and $D_{\mu+k,\nu+l}^{c}$ are isomorphic (\cite{gpots}), the
assumption $\hbar\not = 0$ in Theorem \ref{kth} can be dropped.
\end{rk}

\vspace{.1in}
\begin{thm}
\label{kgroups} $K_{0}(D_{\mu\nu}^{c})\cong Z^{3}+Z_{c}$ and
$K_{1}(D_{\mu\nu}^{c})\cong Z^{3}$.
\end{thm}
\vspace{.1in}
\underline{\em{Proof}}:
\vspace{.1in}
\newline In view of Theorem \ref{kth} and Remark \ref{plusk}, it suffices to
prove the theorem for the commutative case where $D^{c}_{\mu\nu}= C(M_{c})$.
\newline After reparametrizing the Heisenberg group we get that $M_{c}=G/H_{c}$
where
\vspace{.05in}
\[G=\left\{\left( \begin{array}{ccc}
        1 & y & z/c \\
        0 &1 & x  \\
        0 &0 &1
        \end{array}
        \right) : x,y,z\in R \right\}\mbox{\hspace{.1in}and\hspace{.1in}
}H_{c}=\left\{\left( \begin{array}{ccc}
         1 & m & p/c \\
         0 & 1 & q \\
         0 & 0 & 1
         \end{array}
         \right) : m,p,q\in Z \right\}\]
\vspace{.1in}
\newline We first use \cite[Corollary 3]{ro} to reduce the proof to the
computation of the K-theory of $C^{*}(H_{c})$.
\newline The group C*-algebra $C^{*}(H_{c})$ is strong-Morita equivalent to
$C(G/H_{c})\times G$, where G acts by left translation \cite[Example 1]{rfmor}.
Now, G is nilpotent and simply connected so we have
\[G=R\times\!\!\!\mid R\times R\]
as a semi-direct product.

Therefore

\[ C(G/H_{c})\times G \simeq C(G/H_{c})\times R\times\!\!\!\mid R\times R,\]
and Connes'-Thom isomorphism (\cite[10.2.2]{bl}) gives
\[K_{i}(C^{*}(H_{c}))=K_{i}(C(G/H_{c})\times G)
=K_{1-i}(C(G/H_{c}))=K_{1-i}(C(M_{c})).\]

So it suffices to compute $K_{i}(C^{*}(H_{c}))$.
The computation was made  in \cite[Prop. 1.4]{ap} for the case c=1, and the
general case can be obtained with slight modifications to their proof.
 We first write $H_{c}$ as a  semi-direct product, so its group $C^{*}$-algebra
can be expressed as  a crossed product algebra. Then, by using the
Pimsner-Voiculescu exact sequence (\cite[10.2.1]{bl}), we  get its K-groups.
\vspace{.1in}

Let
\[\begin{array}{ccc}
N=\left\{\left( \begin{array}{ccc}
         1 & m & p/c \\
         0 & 1 & 0 \\
         0 & 0 & 1
         \end{array}
         \right) : m,p,\in Z \right\} &
       \mbox{ and }&
K=\left\{\left( \begin{array}{ccc}
         1 & 0 & 0 \\
         0 & 1 & q \\
         0 & 0 & 1
         \end{array}
         \right) : q\in Z \right\} \end{array}\]

Then $H_{c}=N\times_{\alpha_{c}}K$, where $\alpha_{c}$ is conjugation. If we
identify in the obvious way $N$ and $K$ with $Z^{2}$ and $Z$ respectively, we
have that $H_{c}\simeq Z^{2}\times_{\alpha_{c}}Z$, where
$\alpha_{c}(q)(m,p)=(m,p-cmq).$
Then the Pimsner-Voiculescu exact sequence yields:
\[ \begin{array}{ccccc}
   K_{0}(C(T^{2})) & \stackrel{id -\alpha_{c_{*}}}{\longrightarrow} &
K_{0}(C(T^{2})) & \stackrel{i_{*}}{\longrightarrow} & K_{0}(H_{c})    \\
 \begin{array}{cc}
\delta & \uparrow \end{array} &  &  &  &\begin{array}{cc}
\downarrow & \delta \end{array}\\
   K_{1}(H_{c})& \stackrel{i_{*}}{\longleftarrow} & K_{1}(C(T^{2}))&
\stackrel{id-\alpha_{c_{*}}}{\longleftarrow}& K_{1}(C(T^{2}))
\end{array} \]

It was shown on \cite[Prop.1.4]{ap} that $id=\alpha_{1_{*}}$ on
$K_{0}(C(T^{2}))$ and, since $\alpha_{c_{*}}=\alpha_{1_{*}}^{c}$ it follows
that $id=\alpha_{c_{*}}$ on $K_{0}(C(T^{2}))$ for any c. Thus we get the
following short exact sequences:

\[ \begin{array}{ccccccccc}
0 & \longrightarrow & Z^{2} & \longrightarrow & K_{0}(H_{c}) &
\stackrel{\delta}{\longrightarrow} & Ker(id-\alpha_{c_{*}}) & \longrightarrow &
0
\end{array}\]

\[\begin{array}{ccccccccc}
0 & \longrightarrow & K_{1}(C(T^{2}))/Ker(id-\alpha_{c_{*}}) &\longrightarrow &
K_{1}(H_{c})& \stackrel{\delta}{\longrightarrow} &  Z^{2} & \longrightarrow & 0
\end{array},\]
where $id-\alpha_{c_{*}}$ is the map on $K_{1}(C(T^{2}))$.

Let us now compute $id-\alpha_{c_{*}}$ on $K_{1}(C(T^{2}))$.
We have identified $C(T^{2})$ with $C^{*}(Z^{2})$ via Fourier transform, so the
automorphism $\alpha_{c}$ on $C(T^{2})$ becomes $(\alpha_{c}f)(x,y)=f(x-cy,y).$
Now, $K_{1}(C(T^{2}))=Z^{2}$ if we identify $[u_{1}]_{K_{1}}$ and
$[u_{2}]_{K_{2}}$ with (1,0) and (0,1) in $Z^{2}$, respectively, where
$u_{1}(x,y)=e(x)$,\hspace{.1in} $ u_{2}(x,y)=e(y)$\hspace{.1in}for all
$(x,y)\in T^{2}$
Then, for $(a,b)\in Z^{2}$ we have
\[(id-\alpha_{c_{*}})(a,b)=(a,b)-(a,b-ac)=(0,ac).\]
This shows that
\[\begin{array}{cc} Ker(id-\alpha_{c_{*}})=Z\oplus\{0\}\subset Z^{2},&
Im(id-\alpha_{c_{*}})=\{0\}\oplus cZ\subset Z^{2}.
\end{array}\]
So the exact sequences above become:

\[\begin{array}{ccccccccc}
0 & \longrightarrow & Z^{2} & \longrightarrow & K_{0}(H_{c}) & \longrightarrow
& Z & \longrightarrow & 0
\end{array}\]
\[\begin{array}{ccccccccc}
0 & \longrightarrow & Z+Z_{c} & \longrightarrow & K_{1}(H_{c}) &
\longrightarrow &Z^{2} &\longrightarrow & 0.
\end{array}\]
Therefore
\[\begin{array}{ccc} K_{1}(D_{\mu\nu}^{c})=K_{0}(H_{c})=Z^{3}& \mbox{and} &
K_{0}(D_{\mu\nu}^{c})=K_{1}(H_{c})= Z^{3}+Z_{c}.
\end{array}\]
\begin{flushright}
Q.E.D
\end{flushright}
\vspace{.2in}
\underline{\large{\bf Acknowledgement.}} I would like to thank my thesis
adviser, Marc Rieffel, for his constant patience and encouragement, as well as
for a number of helpful suggestions and  comments.
\vspace{.2in}
\begin{flushright}
Centro de Matem\'aticas\\
Facultad de Ciencias\\
Eduardo Acevedo 1139\\
CP 11200\\
Montevideo-URUGUAY.\\
E-mail: abadie@cmat.edu.uy\\
\end{flushright}
\cleardoublepage

\end{document}